\documentclass[10pt,a4paper]{elsarticle}

\usepackage{float}
\usepackage{amstext,amsfonts,amsmath,amssymb}
\usepackage{amsthm}
\usepackage{amsmath}
\usepackage{breqn}
\usepackage{makeidx}
\usepackage{bm}

\makeatletter
\def\ps@pprintTitle{%
  \let\@oddhead\@empty
  \let\@evenhead\@empty
  \def\@oddfoot{\reset@font\hfil\thepage\hfil}
  \let\@evenfoot\@oddfoot
}
\def\@author#1{\g@addto@macro\elsauthors{\normalsize%
\def\baselinestretch{1}%
\upshape\authorsep#1\unskip\textsuperscript{%
\ifx\@fnmark\@empty\else\unskip\sep\@fnmark\let\sep=,\fi
\ifx\@corref\@empty\else\unskip\sep\@corref\let\sep=,\fi}%
\def\authorsep{\unskip,\space}%
\global\let\@fnmark\@empty
\global\let\@corref\@empty 
\global\let\sep\@empty}%
\@eadauthor={#1}
}
\makeatother

\newtheorem{theorem}{Theorem}[section]

\newtheorem{lemma}[theorem]{Lemma}
\newtheorem{corollary}[theorem]{Corollary}

\begin{document}

\renewcommand{\thefootnote}{\fnsymbol{footnote}}

\leftline{\small \tt Preprint submitted to arXiv.org}

\vspace{20.0 mm}

\centerline {\bf \large Accurate estimates of {\boldmath $(1+x)^{1/x}$} Involved in}

\smallskip

\centerline {\bf \large Carleman Inequality and Keller Limit}

\begin{center}
Branko Male\v sevi\' c${}^{\mbox{\tiny $\ast$} \, \mbox{\tiny 1)}}$,
Yue Hu${}^{\mbox{\tiny 2)}}$ and Cristinel Mortici${}^{\mbox{\tiny 3)}}$
\end{center}

\footnotetext{$\!\!\!\!\!\!\!\!\!\!\!{}^{\ast}\,$Corresponding author
}

\footnotetext{$\!\!\!\!\!\!\!\!\!\!\!$E-mails:
Branko~Male\v sevi\' c$\,<${\sl branko.malesevic@etf.bg.ac.rs}$>$,
Yue Hu$\,<${\sl huu3y2@163.com}$>$,
Cristinel~Mortici$\,<${\sl cmortici@valahia.ro}$>$}

\renewcommand{\thefootnote}{\arabic{footnote}}

\begin{center}
{\footnotesize \it
${}^{1)}$ School of Electrical Engineering, University of Belgrade,                                      \\[+0.25 ex]
Bulevar Kralja Aleksandra 73, 11000 Belgrade, Serbia                                                     \\[+1.50 ex]
${}^{2)}$School of Mathematics and Information Science,                                                  \\[+0.25 ex]
Henan Polytechnic University, Jiaozuo, Henan 454000, China                                               \\[+1.50 ex]
${}^{3)}$Valahia University of T\^{a}rgovi\c{s}te, Bd. Unirii 18, 130082 T\^{a}rgovi\c{s}te, Romania;    \\[+0.00 ex]
${}^{3)}$Academy of Romanian Scientists, Splaiul Independen\c{t}ei 54, 050094 Bucharest, Romania;        \\[-0.25 ex]
$\!\!\!\!{}^{3)}$\mbox{University~Politehnica~of~Bucharest,~Splaiul~Independen\c{t}ei~313,~060042~Bucharest,~Romania}
}
\end{center}

\bigskip
\noindent
{\scriptsize
{\bf Abstract.}
In this paper, using the {\rm Maclaurin} series of the functions $(1+x)^{\!1/x}$, some inequalities
from papers \cite{Bicheng_Debnath_1998} and \cite{Mortici_Jang_2015} are generalized.
For arbitrary {\rm Maclaurin} series some general limits of {\rm Keller}'s type are defined and applying
for generalization of some well known results.}

\bigskip
\noindent
{\footnotesize {\bf Keyword.}
Constant $\mathtt{e}$; {\rm Carleman}'s inequality; {\rm Keller}'s limit; approximations

\section{Introduction}

\smallskip Let us start from the function
\begin{equation}
e(x) = \left\{
\begin{array}{ccc}
\left(1+x\right)^{\frac{1}{x}} & : & x \!>\! -1 \,\wedge\, x \!\neq\! 0\,, \\[1.0 ex]
\mathtt{e} & : & x=0\,;%
\end{array}
\right.
\end{equation}
where $\mathtt{e}$ is the base of the natural logarithm. In view of its
importance in study Carleman inequality and Keller's Limit, many interesting
and important results have been obtained for this function and related
functions.

In \cite{Brede_2005}, the following result was established:
\begin{equation}
e(x)=\mathtt{e}\left( 1+\displaystyle\sum_{k=1}^{\infty }{e_{k}\,x^{k}}%
\right) ,
\end{equation}%
for $x\in (-1,1)$ and where sequence $e_{n}$ $\left(n\!\in\!N_{0}\!=\!N\!\cup\!\{0\}\right)$
is determined by
\begin{equation}
\begin{array}{l}
e_{0}=1\,,                                                      \\[0.5 ex]
e_{n}=(-1)^{n}\displaystyle\sum_{k=0}^{n}{\displaystyle\frac{%
(-1)^{n+k}S_{1}(n+k,k)}{(n+k)!}}\displaystyle\sum_{m=0}^{n}{\displaystyle%
\frac{(-1)^{m}}{(m-k)!}}\,;%
\end{array}
\label{Brede}
\end{equation}%
where $S_{1}(p,q)$ is the Stirling's number of the first kind $\left(p, q \!\in\! N\right)$. For the Stirling's numbers of the first kind
the following recurrence relation is true
\begin{equation}
S_{1}(p+1,q) = -\,p S_{1}(p,q) + S_{1}(p,q-1),
\end{equation}
for $p \!>\! q$; $S_{1}(p,p) \!=\! 1$ and $S_{1}(p, q) \!=\! 0$ for $p \!<\! q$
$\left(p, q \!\in\! N\right)$, see for example \cite{Abramowitz_Stegun_1972}, \cite{Mansour_Schork_2015}.
Let us notice that sign of the Stirling's number of the first kind $S_1(p,q)$ is equal to $(-1)^{p-q}$.
Initial values of the sequence of fractions
\begin{equation}
e_{0}\!=\!1,\;
e_{1}\!=\!-\frac{1}{2},\;
e_{2}\!=\!\frac{11}{24},\;
e_{3}\!=\!-\frac{7}{16},\;
e_{4}\!=\!\frac{2447}{5760},\;
e_{5}\!=\!-\frac{959}{2304},\;
e_{6}\!=\!\frac{238043}{580608},\, \ldots
\end{equation}
are tabled with the sequences A055505 and A055535 in \cite{Sloane}. Based on the Lemma 2 from \cite{Brede_2005} (for $t\!=\!1$),
from the representation~(\ref{Brede}) we obtain
\begin{equation}
e_{n}=\mathtt{e}^{-1}\!\displaystyle\sum_{k=1}^{\infty }{\displaystyle\frac{%
S_{1}(n+k,k)}{(n+k)!}},  \label{Brede2}
\end{equation}
where $n \!\in\! N$.


\smallskip
Let us consider the following sequence of partial sums:
\begin{equation}
e_{n}(x) = \mathtt{e}\left(1+\displaystyle\sum_{k=1}^{n}{(-1)^kf_{k}\,x^k}%
\right),
\end{equation}
for
\begin{equation}
f_{k}
=
(-1)^k e_k
=
(-1)^k
\mathtt{e}^{-1}\!
\displaystyle\sum_{i=1}^{\infty }{\displaystyle\frac{S_{1}(k+i,i)}{(k+i)!}};
\end{equation}
where $x \!\in\! (-1,1)$ and $n \!\in\! N_{0}$.

\smallskip
The first aim of the present paper is to establish the monotoneity properties for the $e_{n}(x)$.

\smallskip
In \cite{Mortici_Jang_2015} and \cite {Hu_Mortici_2016}, the authors proved that for any $c \!\in\! R$ the following Keller's type limit holds:
\begin{equation}
\label{Mortici_Lim}
\lim_{n\rightarrow \infty}\left((n+1)\left(1+\frac{1}{n+c}\right)^{n+c}\!-\,n\left(1+\frac{1}{n+c-1}\right)^{n+c-1}\right)=\texttt{e}.
\end{equation}

The second aim of the present paper is to define some general limits of
Keller's type, related to some limits stated in \cite{Mortici_Jang_2015} and
\cite{Hu_Mortici_2016}.

\medskip

\section{Monotoneity properties of the $e_{n}(x)$}

The next statements are true:

\begin{lemma}
The sequence $f_{n}$ is strictly monotonic decreasing.
\end{lemma}

\noindent
\textbf{Proof.}
Based on the recurrence relation for the Stirling's numbers of the first kind the following equalities are true
\begin{equation}
\begin{array}{rcl}
f_{n+1}
\! \! & \!\!=\!\! & \! \!
(-1)^{n+1}\mathtt{e}^{-1}\!\displaystyle\sum_{i=1}^{\infty}{\displaystyle\frac{S_{1}(n+i+1,i)}{(n+i+1)!}} \\[1.0 ex]
\!\! & \!\!=\!\! & \! \!
(-1)^{n+1}\mathtt{e}^{-1}\!\displaystyle\sum_{i=1}^{\infty}{\displaystyle\frac{ -(n+i)S_{1}(n+i,i) + S_{1}(n+i-1,i-1)}{ (n+i+1)!}} \\[1.0 ex]
\!\! & \!\!=\!\! & \!\!
(-1)^{n+1}\mathtt{e}^{-1}\!\displaystyle\sum_{i=1}^{\infty}{\displaystyle\frac{ -(n\!+\!i\!+\!1)S_{1}(n\!+\!i,i) + S_{1}(n\!+\!i,i) + S_{1}(n\!+\!i\!-\!1,i\!-\!1)}{(n+i+1)!}} \\[1.0 ex]
\!\! & \!\!=\!\! & \! \! (-1)^{n}\mathtt{e}^{-1}\!\displaystyle\sum_{i=1}^{\infty}{\displaystyle\frac{S_{1}(n\!+\!i,i)}{(n+i)!}}
+
(-1)^{n+1}\mathtt{e}^{-1}\!\displaystyle\sum_{i=1}^{\infty}{\displaystyle\frac{S_{1}(n\!+\!i,i) + S_{1}(n\!+\!i\!-\!1,i\!-\!1)}{(n+i+1)!}}\,.
\end{array}
\!\!\!\!\!
\end{equation}
Therefore, based on sign of Strirling's numbers of the first kind, we have the conclusion that
\begin{equation}
f_{n} - f_{n+1} = (-1)^{n}\mathtt{e}^{-1}\!\displaystyle\sum_{i=1}^{\infty}{\displaystyle%
\frac{S_{1}(n\!+\!i,i) + S_{1}(n\!+\!i\!-\!1,i\!-\!1)}{(n+i+1)!}} > 0 \, .
\;\;\;\; \Box
\end{equation}
\begin{corollary}
The sequence $f_{n}$ is convergent.
\end{corollary}
\begin{theorem}
Let us form the sequence
\begin{equation}
e_{n}(x) = \mbox{\tt e}\left(1+\displaystyle\sum_{k=1}^{n}{(-1)^kf_{k}\,x^k}%
\right),
\end{equation}
for $x \in (-1,1)$ $\left(n \!\in\! N_{0}\right)$. For $x \!\neq\! 0$ is
true:

\medskip \noindent $(i)$ if $x \in (0,1)$ then
\begin{equation}
e_{1}(x) \!<\! e_{3}(x) \!<\! \ldots \!\!<\! e_{2k-1}(x) \!<\! \ldots
\!\!<\! e(x) \!<\! \ldots \!\!<\! e_{2k}(x) \!<\! \ldots \!\!<\! e_{2}(x)
\!<\! e_{0}(x);
\end{equation}

\noindent $(ii)$ if $x \!\in\! (-1,0)$ then
\begin{equation}
e_{0}(x) \!<\! e_{1}(x) \!<\! e_{2}(x) \!<\! e_{3}(x) \!<\! \ldots \!<\!
e_{2k}(x) \!<\! e_{2k+1}(x) \!<\! \ldots \!<\! e(x).
\end{equation}

\noindent For $x \!=\! 0$ is true that $e_{n}(0) \!=\! e(0) \!=\! \mbox{\tt
e}$ $\left(n \!\in\! N_{0}\right)$.
\end{theorem}

\noindent
\textbf{Proof.}
Theorem obviously follows from the fact that for fixed $x \!\in\! (0,1)$ the following $f_{k}x^k \!>\! 0$ and $f_{k}x^k \!\searrow\! 0$
are~true $\mbox{\small $\Big (\,$} f_{k}\!>\!0, \, f_{k+1}x^{k+1} \!<\! f_{k}x^{k+1} \!= x(f_{k}x^k) \!\leq\! f_{k}x^k \;\mbox{and}\;
\lim\limits_{k \rightarrow \infty}{(f_{k}\,x^{k})} =\! \lim\limits_{k \rightarrow \infty}{f_{k}}
\cdot\!\lim\limits_{k \rightarrow \infty}{x^{k}}\!=\!0\,\mbox{\small $\Big )$}\,.\;\;\Box$
\begin{corollary}
Based on the previous Theorem we obtain proofs of Lemma $2.1$ from {\rm \cite{Bicheng_Debnath_1998}}
and Theorem $1$ from {\rm \cite{Mortici_Jang_2015}} over~$(0,1)$.
\end{corollary}

\section{General limits of Keller's type}

Let us consider an arbitrary Maclaurin series
\begin{equation}
g(x) = \displaystyle\sum_{k=0}^{\infty}{a_{k}x^{k}},
\end{equation}
for $x \!\in\! \left(-\varrho,\varrho\right)$ and some $\varrho \!\in\! \left(0,\infty\right)$,
wherein $a_{k}$ $(k \!\in\! N_{0})$ is some real sequence. Therefore exist a convergent
asymptotic expansion
\begin{equation}
G(y) = g\!\left(\frac{1}{y}\right) = \displaystyle\sum_{k=0}^{\infty}{%
\displaystyle\frac{a_{k}}{y^{k}}}\,,
\end{equation}
for $y \!\in\! R\backslash\{0\}$ and $\mbox{\small $\Big |$}\displaystyle%
\frac{1}{y}\mbox{\small $\Big |$} \!<\! \min\{1,\varrho\}$. Then, the first
general limit of Keller's type of the function~$g$ we define by
\begin{equation}
L_{1} =\! \displaystyle\lim\limits_{y \longrightarrow \infty}{{\Big (}%
(y\!+\!1)G(y)- yG(y\!-\!1){\Big )}}\,.
\end{equation}

We consider the function $(y\!\,+\,\!1)G(y)-yG(y\!\,-\,\!1)$ for values $y
\!>\! 1 \!\,+\,\! \max\left\{1,\frac{1}{\varrho}\right\}$. Using the
binomial expansion $\displaystyle\frac{1}{(y\!-\!1)^k} \!=\! \displaystyle%
\frac{1}{y^k} \left(1-\displaystyle\frac{1}{y}\right)^{\!\!-k}\!\!\!\!  = %
\displaystyle\sum\limits_{i=0}^{\infty}{\displaystyle\frac{{\binom{k+i-1 }{%
k-1}}}{y^{k+i}}}$, for \mbox{$\frac{1}{y} \!<\! 1$} and $k \!\in\! N$, we
obtained, for $y > 1 \!+\! \max\left\{1,\frac{1}{\varrho}\right\}$,
the~following convergent asymptotic expansion
\begin{equation}
\label{Gen-Gen-1}
\begin{array}{rcl}
(y\!+\!1)G(y)- yG(y\!-\!1) \!\! & \!\!=\!\! & \!\! (y\!+\!1) {\bigg (} a_{0}
\!+\! \displaystyle\sum\limits_{k=1}^{\infty}{\displaystyle\frac{a_{k}}{y^{k}%
}} {\bigg )} - y {\bigg (} a_{0} \!+\! \displaystyle\sum\limits_{k=1}^{%
\infty}{\displaystyle\frac{a_{k}}{(y-1)^{k}}} {\bigg )} \\[2.0 ex]
\!\! & \!\!=\!\! & \!\! a_{0} \,+\, \displaystyle\sum\limits_{k=1}^{\infty}{%
\displaystyle\frac{(y\!+\!1)a_{k}}{y^{k}}} \,-\, \displaystyle%
\sum\limits_{k=1}^{\infty}{\displaystyle\sum\limits_{i=0}^{\infty}{\ %
\displaystyle\frac{y{\binom{k+i-1 }{k-1}}a_{k}}{y^{k+i}}}}\,.%
\end{array}%
\end{equation}
For $y \!>\! 1 \!+\! \max\left\{1,\frac{1}{\varrho}\right\}$ let us
determine the representation
\begin{equation}  \label{Gen-Gen-2}
\begin{array}{rcl}
\!\!\!\!\!\!\!\! (y\!+\!1)G(y)\!-\!yG(y\!-\!1) \!\! & \!\!=\!\! & \!\! a_{0}
\!+\! a_{1} \!+\! \displaystyle\frac{a_{1}\!+\!a_{2}}{y} \!+\! \displaystyle
\frac{a_{2}\!+\!a_{3}}{y^2} \!+\! \displaystyle\frac{a_{3}\!+\!a_{4}}{y^3}
\!+\! \displaystyle\frac{a_{4}\!+\!a_{5}}{y^4} \!+\! \ldots               \\[1.75 ex]
\!\! & \!\!-\!\! & \!\! {\Bigg (} \displaystyle\frac{{\binom{0 }{0}}a_{1}}{1}
+ \displaystyle\frac{{\binom{1 }{0}}a_{1}}{y} + \displaystyle\frac{{\binom{2
}{0}}a_{1}}{y^2} + \displaystyle\frac{{\binom{3 }{0}}a_{1}}{y^3} + \ldots \\[1.75 ex]
\!\! & \!\! \!\! & \!\! +\, \displaystyle\frac{{\binom{1 }{1}}a_{2}}{y} + %
\displaystyle\frac{{\binom{2 }{1}}a_{2}}{y^2} + \displaystyle\frac{{\binom{3
}{1}}a_{2}}{y^3} + \displaystyle\frac{{\binom{4 }{1}}a_{2}}{y^4} + \ldots \\[1.75 ex]
\!\! & \!\! \!\! & \!\! +\, \displaystyle\frac{{\binom{2 }{2}}a_{3}}{y^2} + %
\displaystyle\frac{{\binom{3 }{2}}a_{3}}{y^3} + \displaystyle\frac{{\binom{4
}{2}}a_{3}}{y^4} + \displaystyle\frac{{\binom{5 }{2}}a_{3}}{y^5} + \ldots \\[1.75 ex]
\!\! & \!\! \!\! & \!\! +\, \displaystyle\frac{{\binom{3 }{3}}a_{4}}{y^3} + %
\displaystyle\frac{{\binom{4 }{3}}a_{4}}{y^4} + \displaystyle\frac{{\binom{5
}{3}}a_{4}}{y^5} + \displaystyle\frac{{\binom{6 }{3}}a_{4}}{y^6} + \,\ldots
\, \ldots \! {\Bigg )}\,,%
\end{array}
\end{equation}
i.e.
\begin{equation}  \label{Gen-Gen-3}
\begin{array}{rcl}
(y\!+\!1)G(y)- yG(y\!-\!1) \!\! & \!\!=\!\! & \!\! a_{0} - \displaystyle%
\frac{a_1+a_2}{y^2} \\[2.0 ex]
\!\! & \!\! \!\! & \!\!\;\;\;\; - \,\displaystyle\frac{a_1+3a_2+2a_3}{y^3}                   \\[1.50 ex]
\!\! & \!\! \!\! & \!\!\;\;\;\; - \,\displaystyle\frac{a_1+4a_2+6a_3+3a_4}{y^4}              \\[1.50 ex]
\!\! & \!\! \!\! & \!\!\;\;\;\; - \,\displaystyle\frac{a_1+5a_2+10a_3+10a_4+4a_5}{y^5}       \\[1.50 ex]
\!\! & \!\! \!\! & \!\!\;\;\;\; - \,\displaystyle\frac{a_1+6a_2+15a_3+20a_4+15a_5+5a_6}{y^6} \\[-0.75 ex]
\!\! & \!\! \!\! & \!\!\;\;\;\;    \vdots                                                    \\[-0.25 ex]
\!\! & \!\! \!\! & \!\!\;\;\;\; - \,\displaystyle\frac{\;\displaystyle\sum_{i=1}^{k-1}{{
\binom{k }{i - 1}}\,a_{i}} +{\Bigg (}{\binom{k }{k-1}} - 1{\Bigg )}\,a_{k}}{y^k}             \\[-0.25 ex]
\!\! & \!\! \!\! & \!\!\;\;\;\; \vdots%
\end{array}%
\end{equation}
Finally, for $y > 1 + \max\left\{1,\frac{1}{\varrho}\right\}$, we obtained the
following convergent asymptotic expansion
\begin{equation}
\label{Gen-Gen-4}
(y\!+\!1)G(y)- yG(y\!-\!1) = a_{0} - \mathop{\mbox{\Large
$\displaystyle\sum$}}\limits_{k=2}^{\infty}{\ \displaystyle\frac{%
\displaystyle\sum_{i=1}^{k-1}{{\binom{k }{i-1}}\,a_{i}} +{\Bigg (}{\binom{k
}{k-1}} - 1{\Bigg )}\,a_{k}}{y^k}}\,.
\end{equation}
Let us remark that coefficients ${\binom{k }{0}}, {\binom{k }{1}}, {\binom{k}{2}}, ... \,,{\binom{k }{k-1}},
\left(\!{\binom{k }{k-1}} \!-\! 1\!\right)$ respectively to $a_1, a_2, a_3, ...\,,$ $a_{k-1}, a_{k}$ are tabled
as the sequence A193815 in \cite{Sloane}. The above expansion (\ref{Gen-Gen-4}) is sufficient to conclude that
\begin{equation}
L_{1} =\! \displaystyle\lim\limits_{y \longrightarrow \infty}{{\Big (}%
(y\!+\!1)G(y)- yG(y\!-\!1){\Big )}} = a_{0}\,.
\end{equation}
Specially, if $g(x) \!=\! e(x)$, i.e. $G(y) \!=\! g\left(\displaystyle\frac{1}{y}\right) \!=\! \left(1\!+\!\displaystyle\frac{1}{y}\right)^{\!y}$,
then
\begin{equation}
L_{1} =\! \displaystyle\lim\limits_{y \longrightarrow \infty}{\left( {%
\displaystyle\frac{(y\!+\!1)^{y+1}}{y^y} - \displaystyle\frac{y^y}{%
(y\!-\!1)^{y-1}}} \right)} = \mathtt{e}\,.
\end{equation}

\smallskip
Next we use some expansion of type (\ref{Gen-Gen-4}) for generalization of some results for limits of Keller's type
from the paper \cite{Mortici_Jang_2015}. The second general limit of Keller's type of the function~$g$ we define by
\begin{equation}
L_{2} =\! \displaystyle\lim\limits_{y \longrightarrow \infty}{{\Big (}(y\!+\!1)G(y\!+\!c)- yG(y\!+\!c\!-\!1){\Big )}}\,,
\end{equation}
for arbitrary $c \!\in\! R$. Similar to the above discussion, we can obtain

\begin{equation}
\label{Gen_Mortici_Lim}
\begin{array}{l}
(y+1)G(y+c)-yG(y+c-1)                                                                                                   \\[1.0 ex]
=
a_0
+
\mathop{\mbox{\Large $\displaystyle\sum$}}\limits_{k=2}^{\infty}{
\displaystyle\frac{
c \mathop{\mbox{\footnotesize $\displaystyle\sum$}}\limits_{i=1}^{k-1}{\binom{k-1}{i-1}a_i}
-
\left(
\displaystyle\sum_{i=1}^{k-1}{{\binom{k }{i-1}}\,a_{i}} +{\Bigg (}{\binom{k }{k-1}} - 1{\Bigg )}\,a_{k}
\right)}{(y+c)^{k}}}\,.
\end{array}
\end{equation}
From (\ref{Gen_Mortici_Lim}), we can get immediately the limit (\ref{Mortici_Lim}).

\smallskip

\section{Conclusions}

In this paper, the monotoneity properties of the functions $e_{n}(x)$ were
established and some general limits of Keller's type are defined. We believe
that the results will lead to a significant contribution toward the study of
Carleman inequality and Keller's limit.

\break

\bigskip
\noindent
{\bf Acknowledgements.} The first author was supported in part by Serbian Ministry of Education,
Science and Technological Development, Projects ON 174032 and III 44006. The second author was
partially supported by the National Natural Science Foundation of China (no.$\,$11471103).
The third author was partially supported by a Grant of the Romanian National Authority
for Scientific Research, CNCS-UEFISCDI, with the Project Number PN-II-ID-PCE-2011-3-0087.

\bibstyle{plain}

\end{document}